\documentclass[11pt]{article}
\usepackage{amsfonts}
\usepackage{amssymb}
\newtheorem{theorem}{Theorem}[section]

\newtheorem{lemma}[theorem]{Lemma}

\newtheorem{definition}{Definition}[section]
\newtheorem{remark}{Remark}[section]
\newtheorem{example}{Example}[section]
 \usepackage{amssymb}
\usepackage{epsfig}
\usepackage{amsmath}
\usepackage{amsfonts}

\def\tanh#1{\,{\normalsize tanh}{\,#1}\,}

\def\qed{\hbox to 0pt{}\hfill$\rlap{$\sqcap$}\sqcup$\medbreak}

\title{Discontinuous functional differential equations with delayed or advanced arguments}
\author{Rub\'en Figueroa}

\begin{document}

\maketitle

\vspace{1cm}
\begin{center}
{\large Departamento de
An\'alise Matem\'atica,\\
Facultade de Matem\'aticas,\\Universidade de Santiago de Compostela, Campus Sur,\\
15782 Santiago de
Compostela, Spain.
%\\
%Phone: 34 981 56 31 00,
%Ext. 13213 \\ FAX: 34 981 59 70 54
}
\end{center}

\begin{abstract}
We provide new results on the existence of extremal solutions for
discontinuous differential equations with a deviated argument which can be either delayed or advanced. The
boundary condition is allowed to be discontinuous and to depend functionally on the unknown solution.
\end{abstract}

\section{Introduction}

In this paper we discuss the existence of absolutely continuous solutions for the boundary value problem
\begin{equation}\label{p1} \left\{ \begin{array}{ll} x'(t)=f(t,x(t),x(\tau(t)),x) \ \mbox{for almost all (a.a.) }t
\in I=[a,b], \vspace{.2cm}\\

B(x(c),x)=0, \end{array} \right. \end{equation} where $c \in \{a,b\}$ and $\tau:I \longrightarrow I$ is a
measurable function which will satisfy either $\tau(t) \leq t$ for a.a. $ t \in I$, that is, $\tau$ is a {\it
delayed} argument, or $t \leq \tau(t)$ for a.a. $t \in I,$ and in such case we say that $\tau$ is an {\it
advanced} argument. Moreover, functions $f$ and $B$ may be discontinuous in all of their arguments. Finally, we
remark that we will take $c=a$ in delayed problems and $c=b$ in advanced ones. \\

Besides its intrinsic mathematical interest, the convenience of including deviated arguments in mathematical
models of many real--life processes is well--known deal with. A very simple example is the malthusian model of
growth: experience shows that if we replace instant dependence of the form $x'(t)=kx(t)$ by another one of the
form $x'(t)=k(x(t-\tau))$, where $\tau$ is the optimal age of reproduction for the studied species, the model
approaches much better what happens in reality. For a longer discussion about this, we recommend the divulgative
paper \cite{liz}. \\

As start point for our research, we cite the work of Tamasan in \cite{tamasan}. In that paper, the author proved
the existence of extremal solutions for an initial value problem involving an equation with delay, with only
monoto\-ni\-ci\-ty assumptions in the variable involving this deviated argument. In the last years, some authors
looked for new results for this kind of problems, trying to replace monotonicity by less restrictive assumptions
(see for instance \cite{jankowski2} and \cite{jankowski}). In most cases, authors studied functional
di\-ffe\-ren\-tial equations, where deviated problems were only particular cases. This is the case of authors as
Jiang and Wei \cite{jiangwei}, Liz and Nieto \cite{liznieto} or Nieto and Rodr\'iguez-L\'opez \cite{nietorosana1},
\cite{nietorosana2}. Mixed--type equations, that is, with both a delay and an advance argument, were also trashed
out, see \cite{dyki}, \cite{ilea}.\\ The present paper follows the line of \cite{jankowski2} and \cite{jankowski}.
In \cite{jankowski}, Jankowski provides some results for problem (\ref{p1}) with delayed argument and boundary
condition $g(x(a),x(b))=0$. The existence of extremal solutions between given lower and upper solutions is
provided there, with the assumptions that $f$, $g$ and $\tau$ are continuous functions, $f$ satisfies a one-sided
Lipschitz condition in the spatial variables and $g$ is monotone on its first variable and one-sided Lipstichz in
the second one. In \cite{jankowski2} the advanced case is treated following analogous techniques. Now, we modify
those results in order to let functions $f$, $B$ and $\tau$ be discontinuous, using a generalized monotone method that can
also be found in \cite{tamasan}. \\

This paper is organized as follows: in section $2$ we provide a new uniqueness result for discontinuous initial
value problems with deviated arguments under strong Lipschitz conditions. In section $3$ we state two required
ma\-xi\-mum principles proved by Jankowsky and then we show our new results for the boundary problem (\ref{p1})
with one example of application. Finally, in section $4$ we provide a method to construct a lower and an upper
solution for problem (\ref{p1}) in a particular case.

\section{The initial value problem}

In this section, as an auxiliar step in order to achieve our goal, we consider the initial value probem with
delayed argument \begin{equation}\label{ivp} x'(t)=g(t,x(t),x(\tau(t))) \ \mbox{for a.a. } t\in I, \quad x(a)=x_a.
\end{equation}

The main result in this section guarantees the existence of a unique solution for problem (\ref{ivp}) under strong
Lipschitz conditions. \begin{theorem}\label{thivp} Let $\tau:I \longrightarrow I$ be a measurable function such
that $\tau(t) \leq t$ for a.a. $t \in I$ and assume that the following conditions hold: \begin{enumerate}
\item[($H_1$)]{For all $x,y \in \mathbb{R}$ the function $g(\cdot,x,y)$ is measurable.} \item[($H_2$)]{There
exists $\psi \in L^1(I)$ such that for a.a. $t \in I$ and all $x,y \in \mathbb{R} $ we have $|g(t,x,y)| \leq
\psi(t)$.} \item[($H_3$)]{There exist nonnegative functions $L_1,L_2 \in L^1(I)$ such that $$
|g(t,x,y)-g(t,\overline{x},\overline{y})| \leq L_1(t) |x-\overline{x}| + L_2(t) |y-\overline{y}|$$ for a.a. $t \in
I$ and all $x,y \in \mathbb{R}$.} \end{enumerate} Then problem (\ref{ivp}) has a unique absolutely continuous
solution. \end{theorem}

\noindent {\bf Proof.} Consider the operator $A:\mathcal{C}(I) \longrightarrow \mathcal{C}(I)$ defined by $$
Ax(t)=x_a+\int_a^t g(s,x(s),x(\tau(s))) \, ds, \ t \in I,$$ which is well defined by vertue by virtue of
conditions $(H_1)$ and $(H_2)$. \\ It is clear that a fixed point of $A$ is also a solution of problem (\ref{ivp})
and vice versa, so it suffices to show that operator $A$ has a unique fixed point. We will do it by application of
the contractive map theorem. \\

For $x \in \mathcal{C}(I)$ we consider the norm $$ ||x||_* = \max_{t\in I}  e^{-\lambda(t)}|x(t)|, \mbox{ where }
\lambda(t)=\int_a^t (L_1(s)+L_2(s)) \, ds,$$ which makes $\mathcal{C}(I)$ become a Banach space.\\

Let $u,v \in \mathcal{C}(I)$. In view of assumption ($H_3$) we obtain \begin{eqnarray*} ||Au-Av||_* &\leq &\max_{t
\in I}\left\{ e^{-\lambda(t)} \int_a^t |g(s,u(s),u(\tau(s)))-g(s,v(s),v(\tau(s)))| \, ds \right\} \\ &\leq
&\max_{t \in I} \left\{ e^{-\lambda(t)} \int_a^t \left[ L_1(s)|u(s)-v(s)| + L_2(s)|u(\tau(s))-v(\tau(s))| \right]
\, ds \right\} \\ &\leq &||u-v||_* \max_{t \in I} \left\{ e^{-\lambda(t)} \int_a^t e^{\lambda(s)} (L_1(s)+L_2(s))
\, ds \right\} \\ &\leq &||u-v||_* \max_{t \in I} \left\{ e^{-\lambda(t)} \left( e^{\lambda(t)}-1 \right) \right\}
\\ &= &||u-v||_* \max_{t \in I} \left( 1- e^{-\lambda(t)}\right) \leq q ||u-v||_*, \end{eqnarray*} with $q =
\left( 1- e^{-||L_1 + L_2||_{L^1(I)}}\right) < 1$. \\

Then operator $A$ has a unique fixed point, which is the unique solution of problem (\ref{ivp}). \qed

\begin{remark} The previous result hold true if we replace the delayed argument for an advanced one and we
consider a final value problem instead of an initial value one. Indeed, let $\tau:I \longrightarrow I$ be such
that $\tau(t) \ge t$ a.e. and consider the problem \begin{equation}\label{fvp} x'(t)=g(t,x(t),x(\tau(t))) \
\mbox{for a.a. } t\in I, \quad x(b)=x_b. \end{equation} Then $x$ is a solution of problem (\ref{fvp}) if and only
if $y(t)=x(-t)$ is a solution of problem \begin{equation}\label{g} y'(t)=h(t,y(t),y(\hat{\tau}(t))) \ \mbox{for
a.a. } t \in [-b,-a], \quad y(-b)=x_b, \end{equation} where $h(t,y,z)=-g(t,y,z)$ and $\hat{\tau}(t)=-\tau(-t)$,
and now problem (\ref{g}) has the form (\ref{ivp}). \end{remark}

\section{Main results and example}

In order to establish our new results on the existence of extremal solutions for problem (\ref{p1}), we need two
maximum principles proved by Jankowski. The first of them, which can be found in \cite{jankowski}, concerns
problems with delay and the second one, which can be found in \cite{jankowski2}, concerns problems with advance.

\begin{lemma}\label{pmaximo1}{\bf \cite[lemma 2.2]{jankowski}} Let $\tau:I \longrightarrow I$ be a measurable
function such that $\tau(t) \leq t$ a.e. on $I$. Assume that $p \in AC(I)$ and that there exist integrable
functions $K$ and $L$, with $L \geq 0$ a.e. on $I$, satisfying the next inequalities:
\begin{equation}\label{ineq1} \left\{ \begin{array}{ll} p'(t) \leq -K(t) p(t) - L(t) p(\tau(t)) \mbox{ for a.a. }
t \in I, \vspace{.2cm}\\ p(a) \leq 0. \end{array} \right. \end{equation} If \begin{equation}\label{integral1}
\int_a^b L(t) e^{\int_{\tau(t)}^{t} K(s) \, ds} dt \leq 1, \end{equation} then $p \leq 0$ on $I$. \end{lemma}

\begin{lemma}\label{pmaximo2}{\bf \cite[lemma 2.1]{jankowski2}} Let $\tau:I \longrightarrow I$ be a measurable
function such that $t \leq \tau(t)$ a.e. on $I$. Assume that $p \in AC(I)$ and that there exist integrable
functions $K$ and $L$, with $L \geq 0$ a.e. on $I$, satisfying the next inequalities:
\begin{equation}\label{ineq2} \left\{ \begin{array}{ll} p'(t) \geq K(t) p(t) + L(t) p(\tau(t)) \mbox{ for a.a. } t
\in I, \vspace{.2cm}\\ p(b) \leq 0. \end{array} \right. \end{equation} If \begin{equation}\label{integral2}
\int_a^b L(t) e^{\int_{t}^{\tau(t)} K(s) \, ds} dt \leq 1, \end{equation} then $p \leq 0$ on $I$. \end{lemma}

\begin{remark} In the original papers, the previous lemmas are proven with the assumption that function $K$ is
continuous; however the same proofs work in the case that $K$ is only integrable, as we assume. In the same way,
in both lemmas we ask the deviated argument to be only measurable, although in the original version continuity was
required. \end{remark}

The last lemma we need is a version of Bolzano's theorem, whose proof can be found in \cite{fp}.

\begin{lemma}\label{bolzano} {\bf \cite[lemma 2.3]{fp}} Let $a,b \in \mathbb R$, $a \leq b$, and let $h:\mathbb R
\longrightarrow \mathbb R$ be such that $h(a) \leq 0 \leq h(b)$ and \begin{equation} \label{bol} \liminf_{z \to
x^-}h(z) \geq h(x) \geq \limsup_{z \to x^+}h(z) \ \mbox{ for all $x \in [a,b]$.} \end{equation} Then there exist
$c_1,c_2 \in [a,b]$ such that $h(c_1)=0=h(c_2)$ and if $h(c)=0$ for some $c\in [a,b]$ then $c_1 \leq c \leq c_2$,
i.e., $c_1$ and $c_2$ are, respectively, the least and the greatest of the zeros of $h$ in $[a,b]$. \end{lemma}

Before introducing our main result for problem (\ref{p1}) we define the concepts of lower and upper solutions for
that problem. We denote by $AC(I)$ the set of absolutely continuos functions on $I$.

\begin{definition}\label{subsol} We say that $\alpha \in AC(I)$ is a lower solution for problem (\ref{p1}) if the
composition $t \in I \to f(t,\alpha(t),\alpha(\tau(t)),\alpha)$ is measurable and the next inequalities hold:
\begin{equation}\label{lower} \left\{ \begin{array}{ll} \alpha'(t) \leq f(t,\alpha(t),\alpha(\tau(t)), \alpha)
\mbox{ for a.a. } t \in I, \vspace{.2cm}\\ B(\alpha(c),\alpha) \leq 0. \end{array} \right. \end{equation} We
define an upper solution for problem (\ref{p1}) analogously, by reversing the previous inequalities.
\end{definition}

Now we are ready to establish the main results.

\begin{theorem}\label{main1} Let $\tau:I \longrightarrow I$ be a measurable function such that $\tau(t) \leq t$
a.e. Suppose that there exist $\alpha, \beta \in AC(I)$ which are, respectively, a lower and an upper solution for
problem (\ref{p1}) with $\alpha \leq \beta$ on $I,$ and assume that for $f:I \times \mathbb{R}^2 \times AC(I)
\longrightarrow \mathbb{R}$ and $B:\mathbb{R} \times AC(I) \longrightarrow \mathbb{R}$ the following conditions hold: \begin{enumerate}
\item[$(H_1)$]{There exists $\psi \in L^1(I)$ such that for a.a. $t \in I$, all $x \in [\alpha(t),\beta(t)]$, all
$y \in [\alpha(\tau(t)),\beta(\tau(t))]$ and all $\gamma \in [\alpha,\beta]=\{\xi \in AC(I) \, : \, \alpha(t) \leq
\xi(t) \leq \beta(t) \mbox{ on } I \}$ we have $|f(t,x,y,\gamma)| \leq \psi(t)$.} \item[$(H_2)$]{There exist
integrable functions $K, L$, with $L \geq 0$ a.e. on $I$, satisfiying (\ref{integral1}) and such that $$
f(t,x,y,\gamma)-f(t,\overline{x},\overline{y},\overline{\gamma}) \leq K(t)(\overline{x}-x)+L(t)(\overline{y}-y)$$
if $\alpha(t) \leq x \leq \overline{x} \leq \beta(t)$, $\alpha(\tau(t)) \leq y \leq \overline{y} \leq
\beta(\tau(t))$ and $\alpha \leq \gamma \leq \overline{\gamma} \leq \beta$.} \item[$(H_3)$]{For all $\xi \in
[\alpha,\beta]$ and all $x \in \mathbb R$ we have $$ \liminf_{y \to x^-} B(y,\xi) \geq B(x,\xi) \geq \limsup_{y\to
x^+} B(y,\xi),$$ and $B(x,\cdot)$  is nonincreasing in $[\alpha,\beta].$} \end{enumerate} Then problem (\ref{p1})
has extremal solutions in $[\alpha,\beta]$. \end{theorem}

\noindent {\bf Proof.} Let's consider the operator $G:[\alpha,\beta] \longrightarrow AC(I)$ such that maps each
$\xi \in [\alpha,\beta]$ to $G\xi$ defined as the solution of the initial value problem
\begin{equation}\label{ivpaux} \left\{ \begin{array}{ll}
x'(t)=f(t,\xi(t),\xi(\tau(t)),\xi)-K(t)[x(t)-\xi(t)]-L(t)[x(\tau(t))-\xi(\tau(t))], \vspace{.2cm}\\ x(a)=x_{\xi},
\end{array} \right. \end{equation} where $x_{\xi}$ is the greatest solution of the algebraic equation
$B(x,\xi)=0$. \\

\noindent {\it Claim 1: Operator $G$ is well--defined.} Due to hypothesis $(H_3)$ and lemma \ref{bolzano}, the
number $x_{\xi}$ is well--defined. On the other hand, by theorem \ref{ivp} it is clear that $G\xi$ is a
well--defined absolutely continuous function. \\

\noindent {\it Claim 2: $G$ is a nondecreasing operator which maps $[\alpha,\beta]$ into itself.} Let $\xi_1,\xi_2
\in [\alpha,\beta]$ with $\xi_1 \leq \xi_2$. First of all, notice that by condition $(H_3)$ it is $G\xi_1(a) \leq
G\xi_2(a)$. Indeed, by definition it is $G\xi_1(a)=x_{\xi_1}, G\xi_2(a)=x_{\xi_2}$ and by the monotoniciy of $B$
on its second variable we have $$ 0=B(x_{\xi_1},\xi_1) \geq B(x_{\xi_1}, \xi_2) \ \mbox{ and} $$ $$ B(\beta(a),
\xi_2) \geq B(\beta(a), \beta) \geq 0,$$ so by application of lemma \ref{bolzano} function $B(\cdot,\xi_2)$ has at
least one zero in $[x_{\xi_1},\beta]$, and then $x_{\xi_1} \leq x_{\xi_2}$. \\

On the other hand, by condition $(H_2)$ we have $$ (G\xi_1-G\xi_2)'(t) \leq
-K(t)[G\xi_1(t)-G\xi_2(t)]-L(t)[G\xi_1(\tau(t))-G\xi_2(\tau(t))],$$ so by lemma \ref{pmaximo1} we have $G\xi_1
\leq G\xi_2$ on $I$ and then operator $G$ is nondecreasing. The same argument shows that $G\beta \leq \beta$ on
$I$, so $G$ maps $[\alpha,\beta]$ into itself. \\

\noindent {\it Claim 3: Operator $G$ have the extremal fixed points, which are extremal solutions of problem
(\ref{p1}).} We have already proven that $G$ is a nondecreasing operator from the functional interval
$[\alpha,\beta]$ into itself. Moreover, notice that for $\xi \in [\alpha, \beta]$ we have $$ |G\xi'(t)|\leq
\psi(t) + (|K(t)| + L(t))(\beta(t)-\alpha(t)), $$ where function in the right-hand side is integrable on $I$
(notice hat $\beta - \alpha \in L^{\infty}(I)$), so by application of \cite{heikkila}, proposition 1.4.4, we
obtain that $G$ has the extremal fixed points, that is, the greatest, $x^*$, and the least one, $x_*$, with
\begin{equation}\label{maximo} x^*=\max\{x\in [\alpha, \beta] \, : \, x \leq Gx\}, \quad x_*=\min\{x\in [\alpha,
\beta] \, : \, Gx \leq x\}. \end{equation} Let's show that $x_*$ is the least solution of (\ref{p1}) in
$[\alpha,\beta]$ (in a similar way we could prove that $x^*$ is the greatest one). First of all, by virtue of
definition of operator $G$ it is clear that a fixed point of $G$ is also a solution of (\ref{p1}). On the other
hand, let $\zeta \in [\alpha, \beta]$ be another solution of problem (\ref{p1}). So, it is also clear that
$G\zeta=\zeta$ and by the formula (\ref{maximo}) we obtain that $x_* \leq \zeta$. Then, $x_*$ is the least
solution of (\ref{p1}) in $[\alpha, \beta]$. \qed

\begin{remark} Notice that conditon $(H_2)$ in theorem \ref{main1} holds if, in particular, $f$ in nondecreasing
on its spatial variables. \end{remark}

The next result is the analogous of version of theorem \ref{main1} for the advanced case.

\begin{theorem}\label{main2} Let $\tau:I \longrightarrow I$ be a measurable function such that $t \leq \tau(t)$
a.e. Suppose that there exist $\alpha, \beta \in AC(I)$ which are, respectively, a lower and an upper solution for
problem (\ref{p1}) with $\alpha \leq \beta$ on $I,$ and assume that for $f:I \times \mathbb{R}^2 \times AC(I)
\longrightarrow \mathbb{R}$ and $B:\mathbb{R} \times AC(I) \longrightarrow \mathbb{R}$ the following conditions hold: \begin{enumerate}
\item[$(H_1)$]{There exists $\psi \in L^1(I)$ such that for a.a. $t \in I$, all $x \in [\alpha(t),\beta(t)]$, all
$y \in [\alpha(\tau(t)),\beta(\tau(t))]$ and all $\gamma \in [\alpha,\beta]$ we have $|f(t,x,y,\gamma)| \leq
\psi(t)$.} \item[$(H_2)$]{There exist integrable functions $K, L$, with $L \geq 0$ a.e. on $I$, satisfiying
(\ref{integral2}) and such that $$ f(t,x,y,\gamma)-f(t,\overline{x},\overline{y},\overline{\gamma}) \geq
-K(t)(\overline{x}-x)-L(t)(\overline{y}-y)$$ if $\alpha(t) \leq x \leq \overline{x} \leq \beta(t)$,
$\alpha(\tau(t)) \leq y \leq \overline{y} \leq \beta(\tau(t))$ and $\alpha \leq \gamma \leq \overline{\gamma} \leq
\beta$.} \item[$(H_3)$]{For all $\xi \in [\alpha,\beta]$ and all $x \in \mathbb R$ we have $$ \liminf_{y \to x^-}
B(y,\xi) \geq B(x,\xi) \geq \limsup_{y\to x^+} B(y,\xi),$$ and $B(x,\cdot)$  is nonincreasing in
$[\alpha,\beta].$} \end{enumerate} Then, problem (\ref{p1}) has extremal solutions in $[\alpha,\beta]$.
\end{theorem}

\noindent{\bf Proof.} The proof is the same that in theorem \ref{main1}, now using lemma \ref{pmaximo2} instead of
lemma \ref{pmaximo1}. \qed

\begin{remark} Notice that conditon $(H_2)$ in theorem \ref{main2} holds if, in particular, $f$ in nonincreasing
on its spatial variables. \end{remark}

\begin{remark} Notice that neither $f$ nor $B$ need be conti\-nuous in theorems \ref{main1} and \ref{main2}.
\end{remark}

We finish this section by considering an example of a problem that, as far as we know, cannot be studied by any
paper in the literature.

\begin{example} Let's consider the boundary value problem with advance argument \begin{equation}\label{ex1}
\left\{ \begin{array}{ll} x'(t)=\phi(x(t))+\sin t \, x(\sqrt{t})\equiv Fx(t) \ \mbox{for a.a. }t \in I=[0,1],
\vspace{.2cm}\\ x(1)-x(0)=\lambda, \end{array} \right. \end{equation} where $0 < \lambda < 1$ and $$
\phi(x)=\left\{ \begin{array}{ll} \dfrac{1}{2}x-\left(1-\dfrac{1}{n}\right), \ &\mbox{ if } x \in \left[
1-\dfrac{1}{n},1-\dfrac{1}{n+1}\right), \ n=1,2,\ldots \vspace{.2cm}\\ 1, \ &\mbox{ otherwise }. \end{array}
\right. $$ Let $\alpha(t)=0$ and $\beta(t)=t, \ t \in I$. Then, $$ F\alpha(t)=0=\alpha'(t), \ \mbox{a.e. on } I;
\quad B(\alpha(1),\alpha)=-\lambda < 0,$$ and $$ F\beta(t) \leq \dfrac{1}{2} + \sin 1 < 1=\beta'(t), \ \mbox{a.e.
on } I; \quad B(\beta(1),\beta)= 1 - \lambda > 0,$$ so $\alpha$ and $\beta$ are, respectively, a lower and an
upper solution for problem (\ref{ex1}) which moreover satisfy $\alpha \leq \beta$ on $I$. \\ On the other hand,
conditions $(H_1)$ and $(H_2)$ in theorem \ref{main2} hold with $\psi \equiv 1$, $K \equiv \dfrac{1}{2}$ and
$L(t)=\sin t$. Note that in this case we have $$ \int_0^1 L(t) e^{\int_{t}^{\sqrt{t}} K(s) \, ds} dt \leq \sin 1
\, e^{\frac{1}{2}} < 1.$$ As condition $(H_4)$ also hold, we conclude that problem (\ref{ex1}) has the extremal
solutions in the functional interval $[\alpha,\beta]$.

\end{example}

\section{On the existence of lower and upper solutions in particular cases}

The purpose of this section is to provide a method to obtain a lower and an upper solution under certain
assumptions that we will specify later. It is a well--known fact that the construction of lower and upper
solutions for a concrete problem is not a trivial question and, in most cases, it constitutes the {\it real}
problem. \\ The next method is based in the disquisitions of Tamasan in \cite{tamasan}.

\begin{theorem}\label{const} In problem (\ref{p1}) assume that: \begin{enumerate} \item[$(C_1)$]{$\tau$ is a
delayed argument.} \item[$(C_2)$]{For a.a. $t \in I$, all $x \in \mathbb{R}$ and all $\gamma \in AC(I)$ the function
$f(t,x,\cdot,\gamma)$ is nondecreasing.} \item[$(C_3)$]{There exists a nondecreasing linear functional $\phi:AC(I)
\longrightarrow \mathbb{R}$ such that $B(x(a),x)=x(a)-\phi(x)$ for all $x \in AC(I)$.} \item[$(C_4)$]{For a.a. $t \in I$,
all $x,y \in \mathbb{R}$ and all $\gamma \in AC(I)$ we have $$|f(t,x,y,\gamma)| \leq p(t)h(|x|,|y|),$$ where $p \in
L^1(I,\mathbb{R} ^+)$, $h:[0,\infty) \times [0,\infty) \longrightarrow [0,\infty)$ is a nondecreasing function in both of
its arguments and $$ \int_{0}^{\infty} \dfrac{du}{h(u,u)}=\infty. $$} \end{enumerate} A sufficient condition for
the existence of a lower and an upper solution for problem (\ref{p1}) is that there exist nonnegative numbers
$m,n_{\alpha},n_{\beta}$, with $n_i \leq m, \  i=\alpha,\beta$, such that \begin{equation}\label{desiglu} m
-\phi(w) \geq n_i (1- \phi(1)), \ i=\alpha, \beta, \end{equation} where $w$ is the unique solution of the initial
value problem \begin{equation}\label{ivplu} w'(t)=p(t)h(w(t),w(t)), \ t \in I, \quad w(a)=m. \end{equation} In
that case, $$ \begin{array}{ll} \alpha(t)=-w(t)+n_{\alpha}, \\ \beta(t)=w(t)-n_{\beta} \end{array} $$ are
respectively a lower and an upper solution for problem (\ref{p1}), which moreover satisfy $\alpha \leq \beta$ on
$I$. \end{theorem}

\noindent {\bf Proof.} We will show that $\alpha$ is a lower solution for (\ref{p1}). \\ For a.a. $t \in I$ we
have $$ \alpha'(t)=-w'(t)=-p(t)h(w(t),w(t)) \leq -p(t) h(w(t)-n_{\alpha},w(t)-n_{\alpha}) = $$ $$ =-p(t)
h(-\alpha(t),-\alpha(t)) \leq f(t,\alpha(t),\alpha(t),\alpha) \leq f(t,\alpha(t),\alpha(\tau(t)),\alpha).$$ On the
other hand, $$ B(\alpha(a),\alpha)=-w(a)+n_{\alpha} + \phi(w-n_{\alpha}),$$ so condition $B(\alpha(a),\alpha) \leq
0$ is equivalent to $$ w(a)-\phi(w) \geq n_{\alpha} (1-\phi(1)). $$

It can be proven analogously that $\beta$ is an upper solution. The fact that $n_i \leq m, \ i=\alpha, \beta,$
guarantees that $\alpha \leq \beta$ on $I$. \qed

\begin{remark} Theorem \ref{const} hold true if we replace point $a$ by $b$ everywhere and conditions $(C_1)$ and
$(C_2)$ by the next ones: \begin{enumerate} \item[$(C_1)'$] $\tau$ is an advanced argument. \item[$(C_2)'$] For
a.a. $t \in I$, all $x \in \mathbb{R}$ and all $\gamma \in [\alpha,\beta]$ the function $f(t,x,\cdot,\gamma)$ is
nonincreasing. \end{enumerate} \end{remark}

\begin{example} Consider a boundary value problem with delay as \begin{equation}\label{ex2} \left\{
\begin{array}{ll} x'(t)=x(\tau(t))+\tanh [x\left(\frac{1}{2}\right)], \ t \in I=[0,1], \ 0 \leq \tau(t) \leq t
\mbox{ a.e.}, \vspace{.2cm}\\ x(0)=\int_{0}^{1} \int_{0}^{1} k(t,s) x(s) \, ds \, dt, \end{array} \right.
\end{equation} where $k$ is a nonnegative kernel such that $\int_{0}^{1} \int_{0}^{1} k(t,s) \, ds \, dt = K=
\dfrac{1}{8}$ and, as usual, square brackets $[\cdot]$ means integer part. \\

In this case, it is $f(t,x,y)=y+1$ and then condition $(C_4)$ in theorem \ref{const} holds for $p \equiv 1$ and
$h(x,y)=y+1$. So, if we solve the initial value problem (\ref{ivplu}) we obtain $w(t)=(m+1)e^t-1,$ and now
condition (\ref{desiglu}) says that we need $m,n_i$, $i=\alpha,\beta$, such that $m \ge n_i$ and
\begin{equation}\label{desigex} m-\int_0^1 \int_0^1 k(t,s) [(m+1)e^s-1]\, ds\, dt \ge n_i (1-K). \end{equation} As
$$ \int_0^1 \int_0^1 k(t,s) [(m+1)e^s-1] \le [(m+1)e-1] K$$ and $K=\dfrac{1}{8}$, inequality (\ref{desigex}) holds
if $$ m-\dfrac{(m+1)e-1}{8} \ge n_i\left(\dfrac{7}{8}\right).$$ In particular, we can take $m=3$ and $n_i=1$,
$i=\alpha,\beta$, so $$ \alpha(t)=2-4e^t; \ \beta(t)=4e^t-2, \ t \in I,$$ are a lower and an upper solution for
problem (\ref{ex2}) such that $\alpha \leq \beta$ on $I$.

The reader can now check that all conditions in theorem \ref{main1} hold with $\psi(t)=4e^t-1$ and $L\equiv 0$, so
this problem has extremal solutions in $[\alpha,\beta]$.

\end{example}

\begin{center} ACKNOWLEDGEMENT \end{center}
\vspace{.2cm}
The author wants to thank Prof. Rodrigo L. Pouso for his valuable comments and support during the development of
this paper.

%    Bibliographies can be prepared with BibTeX using amsplain, %    amsalpha, or (for "historical" overviews)
%natbib style.
\bibliographystyle{amsplain} %    Insert the bibliography data here.

\end{document}